%common header stuff

\documentstyle{amsppt}
\baselineskip18pt
\magnification=\magstep1
%\NoPageNumbers
%\NoRunningHeads
%\pagewidth{4.5in}
%\pageheight{7.0in}
\pagewidth{30pc}
\pageheight{45pc}

\hyphenation{co-deter-min-ant co-deter-min-ants pa-ra-met-rised
pre-print pro-pa-gat-ing pro-pa-gate
fel-low-ship Cox-et-er dis-trib-ut-ive}
\def\leaderfill{\leaders\hbox to 1em{\hss.\hss}\hfill}

\

\def\idest{i.e.,\ }

\def\a{{\alpha}}
\def\be{{\beta}}
\def\g{{\gamma}}

\def\d{{\delta}}

\def\e{{\varepsilon}}

\def\l{{\lambda}}

\def\w{{\omega}}

\def\ba{{\bold a}}
\def\bb{{\bold b}}

\def\bj{{\bold j}}

\def\bn{{\bold n}}

\def\bu{{\bold u}}
\def\bv{{\bold v}}

\def\b0{\text{\bf 0}}

\def\ra{{\ \longrightarrow \ }}

\def\real{{\Bbb R}}
\def\complex{{\Bbb C}}
\def\zed{{\Bbb Z}}

\def\enn{{\Bbb N}}

\def\End{\text{\rm End}}
\def\Hom{\text{\rm Hom}}

\def\boxit#1{\vbox{\hrule\hbox{\vrule \kern3pt
\vbox{\kern3pt\hbox{#1}\kern3pt}\kern3pt\vrule}\hrule}}
\def\rabbit{\vbox{\hbox{\kern0pt
\vbox{\kern0pt{\hbox{---}}\kern3.5pt}}}}

\def\tableau#1{
        \hbox {
                \hskip -10pt plus0pt minus0pt
                \raise\baselineskip\hbox{
                \offinterlineskip
                \hbox{#1}}
                \hskip0.25em
        }
}

\def\tabCol#1{
\hbox{\vtop{\hrule
\halign{\strut\vrule\hskip0.5em##\hskip0.5em\hfill\vrule\cr\lower0pt
\hbox\bgroup$#1$\egroup \cr}
\hrule
} } \hskip -10.5pt plus0pt minus0pt}

\def\CR{
        $\egroup\cr
        \noalign{\hrule}
        \lower0pt\hbox\bgroup$
}

% Set up the map arrows for commutative diagrams.
%\def\mapright#1{\smash{
%     \mathop{\longrightarrow}\limits^{#1}}}

%Set up macro for commutative diagrams etc. (see Ex. 18.46 in TeXbook)

\def\blank#1#2{
%\hbox to {{#1}}{\vbox to {{#2}}}
\hbox to #1{\hfill \vbox to #2{\vfill}}
}

%Ross's table macros

\def\strut{\vrule height10pt depth5pt width0pt}

\def\fa{{\frak a}}
\def\fb{{\frak b}}
\def\fc{{\frak c}}
\def\fd{{\frak d}}
\def\fe{{\frak e}}
\def\fg{{\frak g}}
\def\fh{{\frak h}}
\def\fgl{{\frak g \frak l}}
\def\seca{1}
\def\secb{2}
\def\secc{3}
\def\secd{4}
\def\sece{5}
\def\secf{6}
\def\secg{7}
\def\sech{8}

\topmatter
\title Representations of Lie algebras arising from polytopes
\endtitle

\author R.M. Green \endauthor
\affil Department of Mathematics \\ University of Colorado \\
Campus Box 395 \\ Boulder, CO  80309-0395 \\ USA \\ {\it  E-mail:}
rmg\@euclid.colorado.edu \\
\newline
\endaffil

\abstract 
We present an extremely elementary construction of the simple Lie algebras
over $\complex$ in all of their minuscule representations, using the vertices
of various polytopes.  The construction itself requires no complicated
combinatorics and essentially no Lie theory other than the definition of a Lie
algebra; in fact, the Lie algebras themselves appear as by-products of
the construction.
\endabstract


\endtopmatter

\centerline{\bf Preliminary version, draft 2}

% check Cooperstein for efficiency savings

\head Introduction \endhead

The simple Lie algebras over the complex numbers are objects of key importance 
in representation theory and mathematical physics.  
These algebras fall into four infinite families ($A_n$, $B_n$, 
$C_n$, $D_n$) and five exceptional types ($E_6$, $E_7$, $E_8$, $F_4$ and
$G_2$).  The classical (\idest non-exceptional) types of Lie algebras are
easily defined in terms of Lie algebras of matrices; such representations
are called the natural representations of the Lie algebras.  However, it is
not so easy to give similar descriptions of the exceptional algebras in a
way that makes it easy to carry out calculations with them.  Another natural
question is whether one can give easy descriptions of other representations
of the classical Lie algebras, such as the spin representations of algebras of
types $B_n$ and $D_n$, which are traditionally constructed in terms of
Clifford algebras (see \cite{{\bf 2}, \S13.5}).

There are several combinatorial approaches to the representation theory
of the simple Lie algebras over $\complex$.  Two of these include Littelmann's
description of representations in terms of paths, and the crystal basis
approach of Kashiwara and the Kyoto school.  Both these approaches are
very versatile but can be combinatorially complicated.  Recent work of 
the author shows how to construct certain Lie algebra
representations using combinatorial structures called ``full heaps'', whose
theory is developed in \cite{{\bf 7}, {\bf 8}}.
The approach of the present paper grew out of an attempt to explain the full
heap representations in as simple a way as possible, and it does not require
any complicated combinatorial constructions.

The polytopes we consider in this paper are convex subsets of $\real^n$
whose vertices (\idest $0$-skeletons) have integer coordinates; such 
polytopes are sometimes called ``lattice polytopes''.  These
include the hypercube, the hyperoctahedron (which is the dual of the
hypercube) and the polytopes known as $2_{21}$ and $3_{21}$ in Coxeter's
notation \cite{{\bf 5}}; the latter two polytopes have $27$ and $56$
vertices respectively.  All these polytopes are highly symmetrical, and the 
symmetry groups have been known for a long time.
The reason that these polytopes are relevant in Lie theory is that
the set of weights for the minuscule representations of simple Lie algebras
over $\complex$ form the vertices of one of the aforementioned polytopes.
This is not obvious, but it is not a complete surprise either: Manivel
\cite{{\bf 13}, Introduction} for example mentions in passing that the weights of 
the $56$-dimensional representation of $\fe_7$ correspond to the vertices of
$3_{21}$.  

Our approach in this paper is to start with the vertices of
the polytope and use them to construct representations of Lie algebras
{\it without} first constructing the Lie algebras themselves.  All the
minuscule representations of simple Lie algebras over $\complex$ may be
constructed in this way, and the construction is remarkably simple.
In \S\secb, we introduce the notion of a ``minuscule system'', which
involves two subsets of $\real^n$, denoted by $\Psi$ and $\Delta$.
The set $\Psi$ is said to be a minuscule
system with respect to the simple system $\Delta$ if two conditions are
satisfied (see Definition \secb.1).  These conditions are very elementary and
easy to check, and whenever they hold, the set $\Delta$ defines a set of linear
operators on a vector space with dimension $|\Psi|$ (Definition \secb.2).
If one makes a judicious choice of $\Psi$ and $\Delta$, then these linear
operators turn out to be the representations of the Chevalley generators
of a simple Lie algebra over $\complex$ acting in one of its minuscule
representations with respect to an obvious basis (the basis can be shown to
be the crystal basis in the sense of \cite{{\bf 12}}, by adapting the
argument of \cite{{\bf 7}, \S8}).
We will show that all minuscule representations 
can be constructed in exactly this way.  The dimension of the space
containing $\Psi$ and $\Delta$ is in some cases much smaller than the
dimension of the representation being constructed and the dimension of the
corresponding Lie algebra.

In all our examples here, $\Psi$
and $\Delta$ are finite sets, and the set $\Delta$ is recognizable as 
either the set of simple roots for a simple Lie algebra, or as the set of 
simple roots together with $\a_0 = -\theta$, where $\theta$ is the highest 
root.  In the latter case, we obtain finite dimensional representations of
certain derived affine Kac--Moody algebras.  Formulating the results in terms 
of affine algebras can be more natural, as the affine algebras have a greater
degree of symmetry.  Another advantage is that it is easier to see how
the modules behave under restriction; for example, the $56$-dimensional
module for the Lie algebra of type $E_7$, after inflation to a module for 
the derived affine algebra, can be restricted to a module for the Lie algebra of
type $A_7$ which is the direct sum of two nonisomorphic $28$-dimensional
irreducible submodules.  Once this observation is made,
our approach here to $\fe_7$ is seen to be very natural.

The layout of this paper is as follows.  In \S\seca, we recall some of the 
basic theory of representations of Lie algebras.  Minuscule systems are
defined in \S\secb, and developed in \S\secc.  Our main result is
Theorem \secc.2, and the linear operators used in it are defined in
Definition \secb.2.  Sections \secd--\secf\  are devoted to examples of
minuscule systems.  \S\secd\  describes the representations of $\fe_6$ and
$\fe_7$ arising from the polytopes $2_{21}$ and $3_{21}$ respectively.  
\S\sece\  describes
various representations arising from the hypercube, including the spin
representations of the Lie algebras of types $B_n$ and $D_n$.  The minuscule
representations of type $A_{n-1}$ are obtained by restricting the spin
representation in type $B_n$ to a subalgebra.  \S\secf\  describes 
representations arising from the hyperoctahedron, including the natural
representations of the Lie algebras of types $C_n$ and $D_n$.  In \S\secg,
we explore connections with algebraic geometry, and concluding remarks
are given in \S\sech.

\head \seca. Background on Lie algebras \endhead

A {\it Lie algebra} is a vector space $\fg$ over a field $k$ equipped with a
bilinear map $[\, ,\ ] : \fg \times \fg \ra \fg$ (the {\it Lie bracket}) 
satisfying the conditions $$\eqalign{
[x, x] &= 0,\cr
[[x, y], z] + [[y, z], x] + [[z, x], y] &= 0,
}$$ for all $x, y, z \in \fg$.  (These conditions are known respectively as 
{\it antisymmetry} and the {\it Jacobi identity}.)

If $\fg_1$ and $\fg_2$ are Lie algebras over a field $k$, then a 
{\it homomorphism} of Lie algebras from $\fg_1$ to $\fg_2$ is a $k$-linear 
map $\phi : \fg_1 \ra \fg_2$ such that $\phi([x, y]) = [\phi(x), \phi(y)]$ for
all $x, y \in \fg_1$.  An {\it isomorphism} of Lie algebras is a bijective
homomorphism.

If $V$ is any vector space over $k$ then the Lie algebra $\fgl(V)$ is the
$k$-vector space of all $k$-linear maps $T : V \ra V$, equipped with the
Lie bracket satisfying $$
[T, U] := T \circ U - U \circ T
,$$ where $\circ$ is composition of maps.

A {\it representation} of a Lie algebra $\fg$ over $k$ is a homomorphism
$\rho : \fg \ra \fgl(V)$ for some $k$-vector space $V$.  In this case,
we call $V$ a (left) {\it module} for the Lie algebra $\fg$ (or a $\fg$-module,
for short) and we say that $V$
{\it affords} $\rho$.  If $x \in \fg$ and $v \in V$, we write $x . v$ to
mean $\rho(x)(v)$.  The {\it dimension} of a module (or of the corresponding
representation) is the dimension of $V$.  If $\rho$ is the zero map, then
the representation $\rho$ and the module $V$ are said to be {\it trivial}.

A {\it submodule} of a $\fg$-module $V$ is a $k$-subspace $W$ of $V$ such
that $x.w \in W$ for all $x \in \fg$ and $w \in W$.  If $V$ has no submodules
other than itself and the zero submodule, then $V$ is said to be 
{\it irreducible}.  

If $V_1$ and $V_2$ are $\fg$-modules, then a $k$-linear map $f : V_1 \ra V_2$
is called a {\it homomorphism of $\fg$-modules} if $f(x.v) = x.f(v)$ for
all $x \in \fg$ and $v \in V_1$.  An {\it isomorphism of $\fg$-modules} is 
an invertible homomorphism of $\fg$-modules.

A subspace $\fh$ of $\fg$ is called a {\it subalgebra} of $\fg$ if
$[\fh, \fh] \subseteq \fh$.  If, furthermore, we have $[\fg, \fh] \subseteq
\fh$ (respectively, $[\fh, \fg] \subseteq \fh$) then $\fh$ is said to be
a {\it left ideal} (respectively, {\it right ideal}) of $\fg$.  If $\fh$ is
both a left ideal and a right ideal of $\fg$, then we call $\fh$ a
{\it two-sided ideal} (or ``ideal'' for short) of $\fg$.  If $\fg$ has no
ideals other than itself and the zero ideal, then $\fg$ is said to be
{\it simple}.  The {\it derived algebra}, $\fg'$, of $\fg$ is the subalgebra 
generated by all elements $\{[x_1, x_2]: x_1, x_2 \in \fg\}$.  It can be
shown that $\fg'$ is an ideal of $\fg$.

\definition{Definition \seca.1}
Let $A$ be an $n \times n$ matrix with integer entries.  Following 
\cite{{\bf 11}, \S1.1}, we call $A = (a_{ij})$ a {\it generalized Cartan matrix} 
if it satisfies the following three properties:
\item{(i)}{$a_{ii} = 2$ for all $1 \leq i \leq n$;}
\item{(ii)}{$a_{ij} \leq 0$ for all $i \ne j$;}
\item{(iii)}{$a_{ij} = 0 \Rightarrow a_{ji} = 0$.}

We call the matrix $A$ {\it symmetrizable} if there exists an invertible
matrix $D$ and a symmetric matrix $B$ such that $A = DB$.
\enddefinition

The next result is a well known presentation for the derived algebra
of a Kac--Moody algebra corresponding to a symmetrizable Cartan matrix.

\vfill\eject

\proclaim{Theorem \seca.2}
Let $A$ be a symmetrizable generalized Cartan matrix.
The derived Kac--Moody algebra $\fg = \fg'(A)$ 
corresponding to $A$ is the Lie algebra over $\complex$ generated
by the elements $\{ e_i, f_i, h_i : i \in \Delta \}$ subject
to the defining relations $$\eqalign{
[h_i, h_j] &= 0,\cr
[h_i, e_j] &= A_{ij} e_j,\cr
[h_i, f_j] &= - A_{ij} f_j,\cr
[e_i, f_j] &= \d_{ij} h_i,\cr
\underbrace
{[e_i, [e_i, \cdots [e_i,} 
_{1 - A_{ij} \text{ times}}
e_j] \cdots ]] &= 0,\cr
\underbrace
{[f_i, [f_i, \cdots [f_i,} 
_{1 - A_{ij} \text{ times}}
f_j] \cdots ]] &= 0,\cr
}$$ where $\d$ is the Kronecker delta.
\endproclaim

\demo{Proof}
This is a special case of \cite{{\bf 11}, Theorem 9.11}.
\qed\enddemo

\remark{Remark \seca.3}
In this paper, we are mostly interested the case where $A$ is of finite
type (as defined in \cite{{\bf 11}, \S4.3}).  In this case, the resulting 
algebra $\fg$ is simple.
\endremark

Suppose for the rest of \S\seca\ that $\fg$ is an algebra satisfying
the hypotheses of Theorem \seca.2.  Let $\fh$ be the subalgebra of $\fg$
spanned by the elements $\{h_i : i \in \Delta\}$.  Let $\fh^* = 
\Hom(\fh, \complex)$ be the dual vector space of $\fh$, and let
$\{\w_i : i \in \Delta\}$ be the basis of $\fh^*$ dual to $\{h_i : i \in
\Delta\}$.
Let $V$ be a $\fg$-module.  An element $v \in V$ is called a 
{\it weight vector} of {\it weight} $\lambda \in \fh^*$ if for all 
$h \in \fh$, we
have $h . v = \l(h) v$.  The weights $\w_i$ are known as {\it fundamental
weights}.
If the weight vector $v$ is annihilated by the
action of all of the elements $e_i$ (respectively, all of the elements $f_i$), 
then we call $v$ a {\it highest weight vector} (respectively, a {\it lowest
weight vector}).

The following result is well known.

\vfill\eject

\proclaim{Proposition \seca.4}
\item{\rm (i)}{Let $\fg$ be a simple Lie algebra over $\complex$.  If
$\l$ is a nonnegative $\zed$-linear combination of the fundamental weights
$\w_i$ then up to isomorphism there is a unique finite dimensional 
irreducible $\fg$-module $L(\l)$ of the form $\fg . v_{\l}$, where $v_{\l}$
is of weight $\l$ and is the unique nonzero highest weight vector of $L(\l)$.  
The modules $L(\l)$ are pairwise nonisomorphic and exhaust all finite 
dimensional irreducible modules of $\fg$.}
\item{\rm (ii)}{Suppose that $V$ is a finite dimensional $\fg$-module 
containing a nonzero highest weight vector $v_\l$ of weight $\l$, and that 
$\dim(V) = \dim(L(\l))$.  Then $V \cong L(\l)$.}
\endproclaim

\demo{Proof}
Part (i) is a special case of \cite{{\bf 2}, Theorem 10.21}.  For part (ii),
it follows from the proof of \cite{{\bf 2}, Proposition 10.13} that any
$\fg$-module $\fg.v_\l$ generated by a highest weight vector of weight
$\l$ is a quotient of the Verma module $M(\l)$.  The Verma module has a unique
maximal submodule, $J(\l)$ (see \cite{{\bf 2}, Theorem 10.9}) and we have
$L(\l) = M(\l)/J(\l)$ by definition.  It follows that $\fg.v_\l$ has 
a quotient module isomorphic to $L(\l)$.  The assumption about dimensions
allows this only if $V \cong L(\l)$ (and $\fg.v_\l = V$).
\qed\enddemo

If $\l$ is a fundamental weight, the corresponding module $L(\l)$ is called
a {\it fundamental module}.  Certain of the fundamental modules for simple Lie
algebras are known as {\it minuscule modules}, for reasons we will not go
into (although see \cite{{\bf 1}, 2.11.15} for an explanation).  The purpose 
of this paper
is to provide a uniform and very elementary construction of these modules.
We now list the minuscule modules, their weights and their dimensions;
more information on this may be found in \cite{{\bf 2}, \S13}.
Our indexing of the weights in this paper is based on that of Kac 
\cite{{\bf 11}}, and in some cases, this differs from Carter's 
notation in \cite{{\bf 2}}.

For the simple Lie algebra of type $A_n$, all the fundamental modules
$$L(\w_1), \ldots, L(\w_n)$$ are minuscule, and we have $$
\dim(L(\w_i)) = {{n+1} \choose {i}}
.$$  In this case, $L(\w_1)$ is the natural module, and $L(\w_i)$ is the
$i$-th exterior power of $L(\w_1)$.

For the simple Lie algebra of type $B_n$ (for $n \geq 2$), 
the only minuscule module is the spin module, $L(\w_n)$, which has 
dimension $2^n$.

For the simple Lie algebra of type $C_n$ (for $n \geq 2$), 
the only minuscule module is the natural module, $L(\w_1)$, which has 
dimension $2n$.

For the simple Lie algebra of type $D_n$ (for $n \geq 4$), 
there are three minuscule modules.  These are the natural module
$L(\w_1)$, of dimension $2n$, and the two spin modules $L(\w_{n-1})$ and
$L(\w_n)$, each of which has dimension $2^{n-1}$.

The simple Lie algebra of type $E_6$ has two minuscule modules, $L(\w_1)$
and $L(\w_5)$, each of which has dimension $27$.

The simple Lie algebra of type $E_7$ has one minuscule module, $L(\w_6)$,
which has dimension 56.

The simple Lie algebras of types $E_8$, $F_4$ and $G_2$ have no minuscule
modules.

\head \secb. Minuscule systems \endhead

\definition{Definition \secb.1}
Let $\Psi$ and $\Delta$ be subsets of vectors in 
$\real^n$ for some $n \in \enn$, where
$\real^n$ is equipped with the usual scalar product and $0 \not\in \Delta$.  We
say that $\Psi$ is a {\it minuscule system with respect to the 
simple system} $\Delta$
if the following conditions are satisfied for every $\bv \in \Psi$ and
$\ba \in \Delta$.
\item{\rm (i)}{We have $2 \bv . \ba = c \ba . \ba$ for some $c = c(\bv, \ba) 
\in \{-1, 0, +1\}$.}
\item{\rm (ii)}{Let $c = c(\bv, \ba)$ be as in (i).  Then we have 
$\bv + \ba \in \Psi$ if
and only if $c = -1$, and we have $\bv - \ba \in \Psi$ if and only if $c = 1$.
(In particular, if $c = 0$, then neither vector $\bv \pm \ba$ lies in $\Psi$.)}
\enddefinition

\definition{Definition \secb.2}
Let $\Psi$ be a minuscule system with respect to the simple system $\Delta$,
and let $k$ be a field.  We define $V_\Psi$ to be the $k$-vector space 
with basis $\{b_\bv : \bv \in \Psi\}$.  For each $\ba \in \Delta$, 
we define $k$-linear endomorphisms
$E_\ba$, $F_\ba$, $H_\ba$ of $V_\Psi$ by specifying their effects on basis 
elements, as follows: $$\eqalign{
E_\ba (b_\bv) &= \cases
b_{\bv + \ba} & \text{if } \bv + \ba \in \Psi;\cr
0 & \text{otherwise};\cr
\endcases\cr
F_\ba (b_\bv) &= \cases
b_{\bv - \ba} & \text{if } \bv - \ba \in \Psi;\cr
0 & \text{otherwise};\cr
\endcases\cr
H_\ba (b_\bv) &= c(\bv, \ba) b_\bv = 2 {{\bv . \ba} \over {\ba . \ba}} 
b_\bv.\cr}$$
\enddefinition

\definition{Definition \secb.3}
Let $\Psi$ be a minuscule system with respect to the simple system $\Delta$.
We define the {\it generalized Cartan matrix}, $A$, of $\Delta$ to be the 
$|\Delta| \times |\Delta|$ matrix whose $(\ba, \bb)$ entry is given by $$
A_{\ba, \bb} = 2 {{\ba . \bb} \over {\ba . \ba}}
.$$
\enddefinition

Although we have apparently given two meanings to the term ``generalized
Cartan matrix'' (the above meaning and Definition \seca.1), they coincide
in all the examples of this paper.  A formulation very similar to Definition
\secb.3 may be found in \cite{{\bf 11}, \S2.3}.

\vfill\eject

\head \secc. Results on minuscule systems \endhead

The following lemma is the key ingredient for our main result.

\proclaim{Lemma \secc.1}
Using the notation of Definition \secb.2, we have the following identities
in $\End(V_\Psi)$, where $\ba, \bb \in \Delta$: 
$$\eqalignno{
H_\ba \circ E_\bb &= E_\bb \circ H_\ba + A_{\ba, \bb} E_\bb, & (1) \cr
H_\ba \circ E_\ba &= E_\ba = -E_\ba \circ H_\ba, & (2) \cr
H_\ba \circ F_\bb &= F_\bb \circ H_\ba - A_{\ba, \bb} F_\bb, & (3) \cr
H_\ba \circ F_\ba &= -F_\ba = -F_\ba \circ H_\ba, & (4) \cr
E_\ba \circ F_\bb &= F_\bb \circ E_\ba = 0 \text{ if } A_{\ba, \bb} < 0, 
&(5) \cr
E_\ba \circ F_\bb &= F_\bb \circ E_\ba \text{ if } A_{\ba, \bb} = 0, 
&(6) \cr
E_\ba \circ E_\ba &= 0, & (7) \cr
E_\ba \circ E_\bb &= E_\bb \circ E_\ba \text{ if } A_{\ba, \bb} = 0, & (8) \cr
E_\ba \circ E_\bb \circ E_\ba &= 0 \text{ if } A_{\ba, \bb} = -1, & (9) \cr
F_\ba \circ F_\ba &= 0, & (10) \cr
F_\ba \circ F_\bb &= F_\bb \circ F_\ba \text{ if } A_{\ba, \bb} = 0, & (11) \cr
F_\ba \circ F_\bb \circ F_\ba &= 0 \text{ if } A_{\ba, \bb} = -1, & (12) \cr
}$$
\endproclaim

\demo{Proof}
We prove (1) by acting each side of the equation on a basis vector $b_\bv$.
If $E_\bb(b_\bv) = 0$, then both sides are trivial, so we may assume this is
not the case, meaning that $\bv + \bb \in \Psi$.  It follows that, in the
notation of Definition \secb.1, we have $c(\bv, \bb) = -1$ and $c(\bv + \bb,
\bb) = 1$.  In turn, this means that $$
H_\ba \circ E_\bb (b_\bv) = H_\ba (b_{\bv + \bb}) 
= 2 {{(\bv + \bb) . \ba} \over {\ba . \ba}} b_{\bv + \bb}
$$ and that $$
E_\bb \circ H_\ba (b_\bv) =
2 {{\bv . \ba} \over {\ba . \ba}} E_\bb(b_\bv)
= 2 {{\bv . \ba} \over {\ba . \ba}} b_{\bv + \bb}
.$$  Subtracting, we have $$
(H_\ba \circ E_\bb - E_\bb \circ H_\ba)(b_\bv) = 
= 2 {{\bb . \ba} \over {\ba . \ba}} b_{\bv + \bb}
= A_{\ba, \bb} E_\bb (b_\bv)
,$$ which proves (1).

If $\bb = \ba$, then the above argument shows that $$
H_\ba \circ E_\ba (b_\bv) = H_\ba (b_{\bv + \ba}) = c(\bv + \ba, \ba) 
b_{\bv + \ba} = E_\ba (b_\bv)
.$$  Part (2) follows from this and the fact that $A_{\ba, \ba} = 2$.

The proof of (3) (respectively, (4)) follows by adapting the argument used
to prove (1) (respectively, (2)).

We now prove that $E_\ba \circ F_\bb = 0$ if $A_{\ba, \bb} < 0$.  By (1)
and (2), we have $$\eqalign{
E_\ba F_\bb &= H_\ba E_\ba F_\bb \cr
&= E_\ba H_\ba F_\bb + 2 E_\ba F_\bb \cr
&= E_\ba F_\bb H_\ba + (2 - A_{\ba, \bb}) E_\ba F_\bb.\cr
}$$  Rearranging, this gives $$
E_\ba F_\bb H_\ba = (A_{\ba, \bb} - 1) E_\ba F_\bb.
$$  Suppose that $E_\ba F_\bb \ne 0$, and let $b_\bv$ be a basis element
for which $E_\ba \circ F_\bb (b_\bv) \ne 0$.  This implies that
$H_\ba(b_\bv) = (A_{\ba, \bb} - 1) b_\bv$, but this is a contradiction to
Definition \secb.1 (i), because $c(\bv, \ba) = A_{\ba, \bb} - 1 \leq -2$.
This shows that $E_\ba \circ F_\bb = 0$, and
the proof that $F_\bb \circ E_\ba = 0$ is very similar, proving (5).

We next turn to (6).  Let us first suppose that $E_\ba \circ F_\bb(b_\bv) 
\ne 0$ for some basis element $\bv$.  (This means that $E_\ba \circ F_\bb
(b_\bv) = b_{\bv - \bb + \ba}$ and that $\bv - \bb + \ba \in \Psi$.) 
By (2) and (3), we have $$\eqalign{
E_\ba \circ F_\bb (b_\bv)
&= - E_\ba \circ H_\ba \circ F_\bb (b_\bv) \cr
&= - E_\ba \circ F_\bb \circ H_\ba (b_\bv). \cr
}$$  It follows that $H_\ba (b_\bv) = -b_\bv$, and that $c(\bv, \ba) = -1$.
In turn, this implies that $\bv + \ba \in \Psi$ and 
$E_\ba (b_\bv) = b_{\bv + \ba}
\ne 0$.  Since $\bv - \bb + \ba \in \Psi$, we have $$
F_\bb \circ E_\ba (b_\bv)
= b_{\bv + \ba - \bb} = E_\ba \circ F_\bb (b_\bv)
.$$  It follows that if $E_\ba \circ F_\bb \ne 0$, then $E_\ba \circ
F_\bb = F_\bb \circ E_\ba$.  The converse statement also follows by a similar
argument.  This in turn implies that $E_\ba \circ F_\bb = 0$ if and only if
$F_\bb \circ E_\ba = 0$, which completes the proof of (6).

The proofs of (8) and (11) follow the same line of argument as the proof of 
(6).

To prove (7), we show that $E_\ba \circ E_\ba (b_\bv) = 0$ for all basis
elements $b_\bv$.  As before, we may reduce to the case where $E_\ba(b_\bv)
\ne 0$, meaning that $\bv + \ba \in \Psi$, $c(\bv, \ba) = -1$ and
$c(\bv + \ba, \ba) = 1$.  The latter fact implies that $E_\ba (b_{\bv + \ba})
= 0$, which completes the proof.  The proof of (10) follows the same argument.

We now prove (9).  As in the proof of (7), the proof reduces to showing
that $$
E_\ba \circ E_\bb \circ E_\ba (b_\bv) = 0
$$ in the case where $c(\bv, \ba) = -1$.  Using (1) and (2), we then
have $$\eqalign{
E_\ba \circ E_\bb \circ E_\ba (b_\bv) 
&= - E_\ba \circ H_\ba \circ E_\bb \circ E_\ba (b_\bv) \cr
&= - E_\ba \circ E_\bb \circ (H_\ba \circ E_\ba (b_\bv))
+ E_\ba \circ E_\bb \circ E_\ba (b_\bv)\cr
&= 0,\cr
}$$ as required.  The proof of (12) follows the same argument as the proof of
(9).
\qed\enddemo

We are now ready to state our main result.

\proclaim{Theorem \secc.2}
Let $\Psi$ be a minuscule system with respect to the simple system $\Delta$,
and let $A$ be the generalized Cartan matrix of $\Delta$.  Assume that $A$ is 
a symmetrizable generalized Cartan matrix in the sense of Definition \seca.1,
and let $\fg$ be the corresponding derived Kac--Moody algebra.  
Then the $\complex$-vector space $V_\Psi$ has the structure of a $\fg$-module,
where $e_i$ (respectively, $f_i$, $h_i$) acts via the endomorphism $E_i$
(respectively, $F_i$, $H_i$).
\endproclaim

\demo{Proof}
We need to show that the defining relations of Theorem \seca.2 are satisfied.

Since the operators $H_\ba$ are simultaneously diagonalizable with respect
to the basis $\{ b_\bv : \bv \in \Psi \}$, they commute, and so we have
$[h_i, h_j] = 0$.

Lemma \secc.1 (1) establishes the relations between the $h_i$ and the $e_j$,
and Lemma \secc.1 (3) establishes the relations between the $h_i$ and the
$f_j$.  Lemma \secc.1 (5) and (6) prove that $[e_i, f_j] = 0$ if $i \ne j$.

We now prove that $[e_i, f_i] = h_i$, for which we need to show that $$
E_i \circ F_i - F_i \circ E_i = H_i
.$$  It is enough to evaluate each side of the equation on a basis element
$b_\bv$.  If $c(\bv, i) = 0$ then all terms act as zero.  If 
$c(\bv, i) = 1$ then $E_i \circ F_i(b_\bv) = b_\bv$, $H_i(b_\bv) = 
b_\bv$, and $F_i \circ E_i(b_\bv) = 0$, thus satisfying the equation.
The case $c(\bv, i) = -1$ is dealt with by a similar argument.

Next we prove that the Serre relation $$
\underbrace
{[e_i, [e_i, \cdots [e_i,} 
_{1 - A_{ij} \text{ times}}
e_j] \cdots ]] = 0
$$ is satisfied.
If $A_{ij} = 0$, this states that $[e_i, e_j] = 0$, which is immediate
from Lemma \secc.1 (8).  If $A_{ij} = -1$, this states that $$
[e_i, [e_i, e_j]] = 0
,$$ or in other words, $$
E_i \circ E_i \circ E_j - E_i \circ E_j \circ E_i + E_j \circ E_i \circ E_i
= 0
,$$ which is immediate from Lemma \secc.1 (7) and (9).  
The only other possibility is that $A_{ij} \leq -2$.  In this
case, every term of the corresponding identity in terms of $E_i$ and $E_j$ 
involves an $E_i \circ E_i$, which is zero by Lemma \secc.1 (7), and this
completes the proof.

A similar argument shows that the Serre relation involving the $f_i$ is 
also satisfied.
\qed\enddemo

The following result provides some methods of constructing new minuscule
systems from known ones, and these will be useful in the sequel.

\proclaim{Proposition \secc.3}
Let $\Psi \subset \real^n$ 
be a minuscule system with respect to the simple system $\Delta$.
Let $\Psi'$ and $\Delta'$ be nonempty subsets of $\Psi$ and $\Delta$,
respectively.
\item{\rm (i)}{Suppose that for every $\bv \in \Psi'$ and $\ba \in \Delta'$,
the following conditions are satisfied.
\item{\rm (a)}{If $c(\bv, \ba) = -1$ then $\bv + \ba \in \Psi'$.}
\item{\rm (b)}{If $c(\bv, \ba) = 1$ then $\bv - \ba \in \Psi'$.}

Then $\Psi'$ is a minuscule system with respect to $\Delta'$.}
\item{\rm (ii)}{If $\Psi' = \Psi$ and $\emptyset \ne \Delta' \subset \Delta$
then $\Psi'$ is a minuscule system with respect to $\Delta'$.}
\item{\rm (iii)}{Let $\bn \in \real^n$ and $l \in \real$.  Suppose that the
sets $$
\Psi(\bn, l) = \{\bv \in \Psi : \bv . \bn = l \}
$$ and $$
\Delta(\bn) = \{\ba \in \Delta : \ba . \bn = 0 \}
$$ are nonempty.
Then $\Psi(\bn, l)$ is a minuscule system with respect to the simple
system $\Delta(\bn)$.}
\endproclaim

\demo{Proof}
Definition \secb.1 applied to $\Psi'$ and $\Delta'$ follows immediately
from the hypotheses of (i).  Part (ii) is an immediate consequence of (i).

Part (iii) follows from (i) and the observation that if $\bv \in \Psi(\bn, l)$
and $\ba \in \Delta(\bn)$ then $(\bv \pm \ba) . \bn = l \pm 0 = l$.
\qed\enddemo

\definition{Definition \secc.4}
If $\Psi'$ and $\Delta'$ satisfy the hypotheses of Proposition \secc.3 (i), 
we will call the pair $(\Psi', \Delta')$ a {\it minuscule subsystem} of 
$(\Psi, \Delta)$.
\enddefinition

We now explain how minuscule systems associated with a Lie algebra also 
support actions of the corresponding Weyl group.

\definition{Definition \secc.5}
Let $n \in \enn$ and $0 \ne \a \in V = \real^n$.  The {\it reflection} $s_\a$
associated to $\a$ is the linear map $s_\a : 
V \ra V$ given by $$
s_\a(\bv) = \bv - 2 {{\bv . \a} \over {\a . \a}} \a
.$$  If $\Psi \subset \real^n$ is a minuscule system with respect to the 
simple system $\Delta$, then we define the Weyl group $W = W_{\Psi, \Delta}$
of $(\Psi, \Delta)$ 
to be the group of automorphisms of $\real^n$ generated by the set
$\{s_\ba : \ba \in \Delta\}$.
\enddefinition

It is not hard to check that this agrees with the usual notion of the Weyl
group associated to a simple Lie algebra over $\complex$ (see 
\cite{{\bf 11}, (1.1.2), \S3.7}).  It is well known \cite{{\bf 10}, \S1.1} that the Weyl 
group action respects the scalar product on $V$.

\proclaim{Proposition \secc.6}
If $\Psi$ is a minuscule system with respect to the 
simple system $\Delta$, then $W = W_{\Psi, \Delta}$ acts on $\Psi$.
\endproclaim 

\demo{Proof}
It is enough to show that if $\ba \in \Delta$ and $\bv \in \Psi$, then
$s_\ba (\bv) \in \Psi$.  By the definitions of $s_\ba$ and $c = c(\bv, \ba)$,
we have $s_\ba(\bv) = \bv - c \ba$, which lies in $\Psi$ by Definition
\secb.1.
\qed\enddemo

\head \S\secd. The Hesse polytope \endhead

In \S\secd, we introduce some examples of minuscule systems related to the
polytope known in Coxeter's notation as $3_{21}$.  
This polytope does not have a consistent name in the literature; we will
follow Conway and Sloane in calling $3_{21}$ the {\it Hesse polytope}, as
this name does not appear to have any other connotations.
The Hesse polytope has $56$ vertices, whose coordinates are given by the
set $\Psi^{E_7}$ of Definition \secd.1. 
Note that we have multiplied Conway and Sloane's coordinates for the vertices
by $4$, in order to make them integers and to retain compatibility with
du Val's coordinates \cite{{\bf 6}, \S7}.

The {\it Schl\"afli polytope}, which is called $2_{21}$ in Coxeter's notation,
also plays a role in the examples of this section involving the Lie algebra
of type $E_6$.  It has $27$ vertices, whose coordinates can be given by
either of the sets $\Psi(\bn, \pm 8)$ appearing in Proposition \secd.3.
More details on the inclusion of the Schl\"afli polytope in the Hesse 
polytope may be found in \cite{{\bf 3}, \S9}.

\definition{Definition \secd.1}
Let $\e_0, \e_1, \ldots, \e_7 \in \real^8$ be such that $\e_i$ has a $1$ in
position $i+1$, and zeros elsewhere.
For $0 \leq i, j \leq 7$, define the vector $\bv_{i, j} = \bv_{\{i, j\}} = 
\bv_{j, i} \in \real^8$ by $$
\bv_{i, j} := 4(\e_i + \e_j) - \left( \sum_{i = 0}^7 \e_i \right)
.$$  (For example, we have $\bv_{0, 1} = (3, 3, -1, -1, -1, -1, -1, -1)$.)
Let $\Psi^{E_7}$ consist of the $56$ vectors $\{\pm \bv_{i, j} : 
1 \leq i < j \leq 8\}$.  

It is convenient for later purposes to introduce the sets 
$K_0 = \{0, 1, 2, 3\}$ and $K_7 = \{4, 5, 6, 7\}$.
\enddefinition

\proclaim{Lemma \secd.2}
Let $\Psi^{E_7}$ be as in Definition \secd.1, and let
$$\Delta^{E_7^{(1)}} = 
\{\a_0, \a_1, \ldots, \a_7\},$$ where $
\a_i = 4(\e_i - \e_{i+1})
$ if $0 \leq i < 7$, and $
\a_7 = (-2, -2, -2, -2, 2, 2, 2, 2)
.$  Then $\Psi^{E_7}$ is a minuscule system with respect to the simple
system $\Delta^{E_7^{(1)}}$.  
\endproclaim

\demo{Proof}
Suppose first that $\ba = \a_i$ for some $i < 7$, and let $\bv \in \Psi^{E_7}$.
Write $\bv = \sum_{j = 0}^7 \l_j \e_j$.
The proof is a case by case check according to the values of $\l_i$ and
$\l_{i+1}$.  There are three cases to check.

The first possibility 
is that $\l_i = \l_{i+1}$.  This implies that $\bv . \a_i = 0$.  The coefficients of $\e_i$ and of $\e_{i+1}$ in $\bv + \ba$ differ by $8$, 
which means that $\bv + \ba \not\in \Psi$, and a similar argument shows that
$\bv - \ba \not\in \Psi$.  The conditions of Definition \secb.1 are therefore
satisfied.

The second possibility 
is that $(\l_i, \l_{i+1}) \in \{(-3, 1), (-1, 3)\}$, that is,
$\l_{i+1} = \l_i + 4$.  This implies that, 
$\bv . \ba = -16$ and $\ba . \ba = 32$.
This satisfies Definition \secb.1 (i) with $c = -1$.  In this case,
$\bv - \ba \not\in \Psi$, because the coefficients of $\e_i$ and $\e_{i+1}$
in $\bv - \ba$ do not lie in the set $\{\pm 3, \pm 1\}$.  However, the
vector $\bv + \ba$ is obtained from $\bv$ by exchanging the coefficients of
$\e_i$ and $\e_{i+1}$, which means that $\bv + \ba \in \Psi$.  This satisfies
Definition \secb.1 (ii).

The third possibility is that $(\l_i, \l_{i+1}) \in \{(3, -1), (1, -3)\}$, 
that is,
$\l_{i+1} = \l_i - 4$.  An analysis like that of the previous paragraph
shows that $c = 1$, $\bv + \ba \not\in \Psi$, and $\bv - \ba \in \Psi$, as
required.

It remains to show that Definition \secb.1 is satisfied with $\ba = \a_7$.
To check this, we use the sets $K_0, K_7$ of Definition \secd.1.  
Let $\bv = \pm \bv_{i, j}$.  As before, there are three cases to check.

The first possibility is that
$\{i, j\} \not\subseteq K_l$ for some $l \in \{0, 1\}$.  (Informally, this 
means that the two occurrences of $3$ (or $-3$) in $\bv$ do not occur in the
same half of the vector.)  This implies that $\bv . \a_7 = 0$.  Furthermore,
neither of the two vectors $\bv \pm \a_7$ lies in $\Psi$, because in 
each of them, one of the basis vectors $\e_i$ appears with coefficient $\pm 5$.
Definition \secb.1 is therefore satisfied in this case.

The second possibility 
is that either $\bv = +\bv_{i, j}$ with $\{i, j\} \subset K_0$,
or that $\bv = -\bv_{i, j}$ with $\{i, j\} \subset K_7$.  In each case,
$\bv . \a_7 = -16$ and $\bv + \a_7 \in \Psi$.  However, in each case, we have
$\bv - \a_7 \not\in \Psi$, because two basis vectors appear in $\bv - \a_7$ 
with coefficient $\pm 5$.  Since $\a_7 . \a_7 = 32$, Definition \secb.1
is satisfied with $c = -1$.

The third possibility 
is that either $\bv = +\bv_{i, j}$ with $\{i, j\} \subset K_7$,
or that $\bv = -\bv_{i, j}$ with $\{i, j\} \subset K_0$.  An analysis like
that of the above paragraph shows that Definition \secb.1 is satisfied
with $c = 1$, $\bv - \a_7 \in \Psi$, and $\bv + \a_7 \not\in \Psi$.  This
completes the proof.
\qed\enddemo

\proclaim{Proposition \secd.3}
Let $\Psi = \Psi^{E_7}$ be as in Definition \secd.1, and let $\Delta = 
\Delta^{E_7^{(1)}}$ be as in Definition \secd.2.
\item{\rm (i)}{The $56$-dimensional $\complex$-vector space $V_\Psi$ 
has the structure of a $\fg$-module, where $\fg$ is the derived affine 
Kac--Moody algebra of type $E_7^{(1)}$.}
\item{\rm (ii)}{Let $\Psi' = \Psi$ and $\Delta' = \Delta \backslash \{\a_0\}$.
Then $\Psi'$ is a minuscule system with respect to the simple system
$\Delta'$, and $V_{\Psi'}$ is a module for the simple Lie algebra $\fe_7$ 
over $\complex$ of 
type $E_7$.  It is an irreducible module with highest weight vector 
$-\bv_{0, 7}$ and lowest weight vector $\bv_{0, 7}$ (as in Definition
\secd.1).}
\item{\rm (iii)}{Let $\bn = \bv_{0, 7}$.  Then we have a disjoint union $$
\Psi = \Psi(\bn, 24) \ 
\dot\cup \ \Psi(\bn, 8) \ 
\dot\cup \ \Psi(\bn, -8) \ 
\dot\cup \ \Psi(\bn, 24) 
.$$  For $l \in \{24, 8, -8, -24\}$,
$\Psi(\bn, l)$ is a minuscule system with respect to the
simple system $\Delta(\bn) = \Delta \backslash \{\a_0, \a_6\}$, and 
$V_{\Psi(\bn, l)}$ is a module for the simple Lie algebra $\fe_6$ over
$\complex$ of type $E_6$.  The two modules $V_{\Psi(\bn, \pm 24)}$ are 
trivial one-dimensional modules for $\fe_6$, whereas the two modules
$V_{\Psi(\bn, \pm 8)}$ are nonisomorphic $27$-dimensional irreducible
modules.  The module $V_{\Psi(\bn, 8)}$ has highest weight
$\bv_{1, 7}$ and lowest weight $\bv_{0, 6}$.  The module $V_{\Psi(\bn, -8)}$
has highest weight $-\bv_{0, 6}$ and lowest weight $-\bv_{1, 7}$.}
\item{\rm (iv)}{For $l \in \{24, 8, -8, -24\}$,
$\Psi(\bn, l)$ is a minuscule system with respect to the
simple system $\Delta(\bn) \cup \{\a\}$, where $\a = 4(\e_7 - \e_0)$.
This makes $V_{\Psi(\bn, l)}$ into a module for the derived affine Kac--Moody
algebra $\fg$ of type $E_6^{(1)}$.}
\endproclaim

\demo{Proof}
By Lemma \secd.2, $\Psi$ is a minuscule system with respect to
the simple system $\Delta$.  One may check directly (using \cite{{\bf 11},
\S2.3}) that the
associated matrix $A$ is the symmetrizable generalized Cartan matrix of
type $E_7^{(1)}$ of \cite{{\bf 11}}.  Theorem
\secc.2 then establishes (i).

For (ii), we know that $\Psi'$ is a minuscule system with respect to the
simple system $\Delta'$ by Proposition \secc.3 (ii).  The matrix $A$ in this
case is the symmetrizable (generalized) Cartan matrix of type $E_7$ of
\cite{{\bf 11}}.  It follows from Theorem \secc.2 that $V_{\Psi'}$ is a
module for $\fe_7$.  
Direct checks show that $-\bv_{0, 7}$ is annihilated
by all the operators $E_i$, $\bv_{0, 7}$ is annihilated by all
the operators $F_i$, and $-\bv_{0, 7}$ is annihilated by all the operators
$H_i$ except $H_6$, in which case we have $H_6(-\bv_{0, 7}) = -\bv_{0, 7}$.
Since $V_{\Psi'}$ has the same dimension as $L(\w_6)$ and contains a highest
weight vector of weight $\w_6$, the modules $V_{\Psi'}$ and $L(\w_6)$ 
are isomorphic and irreducible by Proposition \seca.4.

We next establish the decomposition of $\Psi$ described in (iii).  We
have $\Psi(\bn, 24) = \{\bv_{0, 7}\}$ and $\Psi(\bn, -24) = \{-\bv_{0, 7}\}$.
The set $\Psi(\bn, 8)$ consists of the vectors $$
\{ \bv_{0, i} : 1 \leq i \leq 6 \}
\cup 
\{ \bv_{i, 7} : 1 \leq i \leq 6 \}
\cup
\{ -\bv_{i, j} : 1 \leq i < j \leq 6 \}
,$$ and we have $\Psi(\bn, -8) = - \Psi(\bn, 8)$.  It is easy to check that
$\Psi$ is the disjoint union of these four sets.  Proposition \secc.3 (iii)
shows that $\Psi(\bn, l)$ is a minuscule system with respect to $\Delta(\bn)$,
and Theorem \secc.2 shows that the modules $V_{\Psi(\bn, l)}$ are modules
for $\fe_6$ (after the generalized Cartan matrix has been recognized as
symmetrizable of type $E_6$).  The assertions about dimensions and weight
vectors are easy to check.

A quick calculation shows that, for $$
H_i . \bv_{1, 7} = \cases 
\bv_{1, 7} & \text{if } i = 1,\cr
0 & \text{if } i \in \{2, 3, 4, 5, 7\};\cr
\endcases$$ in contrast, we have $$
H_i . (-\bv_{0, 6}) = \cases
-\bv_{0, 6} & \text{if } i = 5,\cr
0 & \text{if } i \in \{1, 2, 3, 4, 7\}.\cr
\endcases$$  This shows that $V_{\Psi(\bn, 8)}$ (respectively,
$V_{\Psi(\bn, -8)}$ has the same dimension as, and a nonzero weight vector
of the same weight as $L(\w_1)$ (respectively, $L(\w_5)$).  Proposition 
\seca.4 now shows that the two modules $V_{\Psi(\bn, \pm 8)}$
are irreducible and nonisomorphic.

To prove (iv), we need to check that Definition \secb.1 is satisfied with
$\ba = \a$.  This follows by imitating the case
analysis for the case $i < 7$ in Lemma \secd.2, using the fact
that $\bn . \a = 0$.
\qed\enddemo

\head \sece. The hypercube \endhead

In \S\sece, we consider examples relating to the polytope known as the
the {\it hypercube} or {\it measure polytope}; in Coxeter's notation
it is denoted $\g_n$.  The set $\Psi$ defined in Lemma \sece.1 is
our standard set of coordinates for the $2^n$ vertices of the hypercube.

We will show how the hypercube may be used to construct the spin
representations of the simple Lie algebras of types $B_n$ and $D_n$.  By
passing to appropriate subsystems, we obtain all the fundamental 
representations of the simple Lie algebra of type $A_n$ as a by-product.

\proclaim{Lemma \sece.1}
Let $n \geq 3$, let $\e_0, \ldots, \e_{n-1} \in \real^n$ be the usual basis for
$\real^n$, and let $\Psi$ be the set of $2^n$ vectors of the form $$
(\pm 2, \pm 2, \ldots, \pm 2)
.$$  Let $\Delta = \{\a_0, \a_1, \ldots, \a_n\}$, where 
$\a_0 = -4(\e_0 + \e_1)$, $\a_n = 4 \e_{n-1}$, and $\a_i = 4(\e_{i-1}
- \e_i)$ for $0 < i < n$.
Then $\Psi$ is a minuscule system with respect to the simple
system $\Delta$.  
\endproclaim

\demo{Proof}
We check that Definition \secb.1 holds for each of the $\a_i$ in turn.
Let $\bv = \sum_{j = 0}^{n-1} \l_j \e_j \in \Psi$.

Suppose first that $0 < i < n$.
The proof is a case by case check according to the values of $\l_i$ and
$\l_{i+1}$.  There are three cases to check, and we omit the details because
the cases are almost identical to those in 
the first part of the argument proving Lemma \secd.2.

Next, suppose that $i = 0$.  There are three cases to check, according
to the values of $\l_0$ and $\l_1$.  
If $\l_0 = \l_1 = +2$, then $\bv + \a_0 \in \Psi$, $\bv - \a_0 \not\in \Psi$,
and $2 \bv . \a_0 = -32 = -\a_0 . \a_0$, giving $c = c(\bv, \a_0) = -1$ 
as required.
If $\l_0 = \l_1 = -2$, then $\bv - \a_0 \in \Psi$, $\bv + \a_0 \not\in \Psi$,
and $2 \bv . \a_0 = 32 = \a_0 . \a_0$, giving $c = 1$ as required.
If $\l_0 \ne \l_1$, then neither vector $\bv \pm \a_0$ lies in $\Psi$,
and $2 \bv . \a_0 = 0$, giving $c = 0$.  Definition \secb.1 is therefore
satisfied in all three cases.

Finally, suppose that $i = n$.  There are two cases to check, according to
the value of $\l_{n-1}$.
If $\l_{n-1} = +2$ then $\bv - \a_n \in \Psi$ and $\bv + \a_n \not\in \Psi$.
We also have $2 \bv . \a_n = 16 = \a_n . \a_n$, giving $c = 
c(\bv, \a_n) = 1$, thus satisfying Definition \secb.1.
If $\l_{n-1} = -2$ then $\bv + \a_n \in \Psi$ and $\bv - \a_n \not\in \Psi$.
We also have $2 \bv . \a_n = -16 = -\a_n . \a_n$, giving $c = -1$, thus 
satisfying Definition \secb.1 and completing the proof.
\qed\enddemo

We may now state an analogue of Proposition \secd.3.

\proclaim{Proposition \sece.2}
Maintain the notation of Definition \sece.1.
Let $\bj = \sum_{j = 0}^{n-1} \e_j$ and $$
S = \{2n-4j : 0 \leq j \leq n \}
.$$
\item{\rm (i)}{The $2^n$-dimensional $\complex$-vector space $V_\Psi$ 
has the structure of a $\fg$-module, where $\fg$ is the derived affine 
Kac--Moody algebra of type $B_n^{(1)}$.}
\item{\rm (ii)}{Let $\Psi' = \Psi$ and $\Delta' = \Delta \backslash \{\a_0\}$.
Then $\Psi'$ is a minuscule system with respect to the simple system
$\Delta'$, and $V_{\Psi'}$ is a module for the simple Lie algebra $\fb_n$ 
over $\complex$ of 
type $B_n$.  It is an irreducible module with highest weight vector $2\bj$
and lowest weight vector $-2\bj$,
and affords the spin representation of $\fb_n$.}
\item{\rm (iii)}{We have a disjoint union $$
\Psi = \bigcup_{j = 0}^n \Psi(\bj, 2n-4j)
.$$  For $l \in S$,
$\Psi(\bj, l)$ is a minuscule system with respect to the
simple system $$\Delta(\bj) = \Delta \backslash \{\a_0, \a_n\},$$ and 
$V_{\Psi(\bj, l)}$ is a module for the simple Lie algebra $\fa_{n-1}$ over
$\complex$ of type $A_{n-1}$.  The two modules $V_{\Psi(\bj, \pm 2n)}$ are 
trivial one-dimensional modules for $\fa_{n-1}$, and the other modules
$V_{\Psi(\bj, l)}$ satisfy $$
V_{\Psi(\bj, 2n-4j)} \cong L(\w_{n-j})
.$$  The module $V_{\Psi(\bj, 2n-4j)}$ has highest weight $$
- 2 \bj + 4 \left(\sum_{i = 0} ^ {n-j-1} \e_i \right)
$$ and lowest weight $$
2 \bj - 4 \left(\sum_{i = 0} ^ {j-1} \e_i \right)
.$$}
\item{\rm (iv)}{For $l \in S$, $\Psi(\bj, l)$ is a minuscule system with
respect to the simple system $\Delta(\bj) \cup \{\a\}$, where
$\a = 4(\e_{n-1} - \e_0)$.
This makes $V_{\Psi(\bj, l)}$ into a module for the derived affine Kac--Moody
algebra $\fg$ of type $A_{n-1}^{(1)}$.}
\endproclaim

\demo{Proof}
Using Lemma \sece.1 in place of Lemma \secd.2, the proof of (i)
follows the same argument as the proof of Proposition \secd.3 (i).

The proof of (ii) now follows by copying the argument of 
Proposition \secd.3 (ii). 
In this case, the module turns out to be $L(\w_n)$.

It is easily checked that $\Psi(\bj, 2n-4j)$ consists precisely of the
vectors in $\Psi$ that have $j$ occurrences of $-2$, from which the first
assertion of (iii) follows.  Proposition \secc.3 (iii)
shows that $\Psi(\bj, 2n-4j)$ is a minuscule system with respect to 
$\Delta(\bj)$,
and Theorem \secc.2 shows that the modules $V_{\Psi(\bj, 2n-4j)}$ are modules
for $\fa_{n-1}$ (after the generalized Cartan matrix has been recognized as
symmetrizable of type $A_{n-1}$).  The assertions about dimensions and weight
vectors are easy to check.  If $j \ne \pm n$ and $\bv$ is the
highest weight vector $\bv$ of $V_{\Psi(\bj, 2n-4j)}$, then we have
$H_i . \bv = 0$ unless $i = n-j$, in which case $H_i . \bv = \bv$.
The required isomorphism now follows from Proposition \seca.4.

To prove (iv), we may copy the argument of Proposition \secd.3 (iv) to 
check that Definition \secb.1 is satisfied with
$\ba = \a$.  (Note that $\bj . \a = 0$.)
\qed\enddemo

\proclaim{Lemma \sece.3}
Let $n \geq 4$, let 
$\e_0, \ldots, \e_{n-1} \in \real^n$ be the usual basis for
$\real^n$, and let $\Psi$ be as in Lemma \sece.1.
Let $\Psi_D^+$ (respectively, $\Psi_D^-$) be the subset of $\Psi$ whose 
vectors contain an even (respectively, odd) number of occurrences of $-2$.

Let $\Delta_D = \{\a_0, \a_1, \ldots, \a_{n-1}, \a'_n\}$, where 
$\a'_n = 4 (\e_{n-2} + \e_{n-1})$ and the other vectors $\a_i$ are 
as in Lemma \sece.1.

Then $\Psi = \Psi_D^+ \ \dot\cup \ \Psi_D^-$ 
is a minuscule system with respect to the simple
system $\Delta_D$, and both $(\Psi_D^+, \Delta_D)$ and $(\Psi_D^-, \Delta_D)$
are minuscule subsystems of $(\Psi, \Delta_D)$.
\endproclaim

\demo{Proof}
Most of the work for checking that $\Psi$ is a minuscule system with
respect to $\Delta_D$ is done in the proof of Lemma \sece.1.  The
only extra criterion to check is that Definition \secb.1 holds for
$\ba = \a'_n$.  This follows by making appropriate sign changes to the
argument used to check Definition \secb.1 for $\ba = \a_0$ as in the proof
of Lemma \sece.1.

Letting $\bj$ be as in Proposition \sece.2, we see
that $\bv \in \Psi$ lies in $\Psi_D^+$ if the integer $\bv . \bj$ is a 
multiple of $8$, and $\bv$ lies in $\Psi_D^-$ otherwise.  We observe that
each $\ba \in \Delta_D$ has the property that $\ba . \bj$ is a multiple of
$8$.  We now apply Proposition \secc.3 (i), which proves that 
$(\Psi_D^\pm, \Delta_D)$ are minuscule subsystems.
\qed\enddemo

\proclaim{Proposition \sece.4}
Maintain the notation of \sece.1--\sece.3.
\item{\rm (i)}{Each of the $2^{n-1}$-dimensional $\complex$-vector spaces 
$V_{\Psi_D^\pm}$ 
has the structure of a $\fg$-module, where $\fg$ is the derived affine 
Kac--Moody algebra of type $D_n^{(1)}$.}
\item{\rm (ii)}{Let $\Psi^\pm = \Psi_D^\pm$ and 
$\Delta^\pm = \Delta_D \backslash \{\a_0\}$.
Then each of the two sets $\Psi^\pm$ is a minuscule system with 
respect to each of the simple systems $\Delta^\pm$ respectively, 
and each of the two spaces $V_{\Psi^\pm}$ is a module for the simple Lie 
algebra $\fd_n$ over $\complex$ of 
type $D_n$.  The modules are nonisomorphic and both irreducible, and they 
afford the two spin representations of $\fd_n$.  
The highest weight vector of $V_{\Psi^+}$ (respectively, $V_{\Psi^-}$) is 
$2\bj$ (respectively, $2\bj - 4\e_{n-1}$).
The lowest weight vectors of $V_{\Psi^\pm}$ are $-2\bj$ and
$-2\bj + 4\e_{n-1}$, where the assignment of vectors to modules depends 
on whether $n$ is even or odd.}
\endproclaim

\demo{Proof}
Using Lemma \sece.3 in place of Lemma \secd.2, the proof
of (i) follows the same argument as the proof of Proposition \secd.3 (i).

The first assertion of (ii) follows by using Lemma \sece.3 and
copying the argument of Proposition \secd.3 (ii).  
The operators $H_i$ (for $1 \leq i < n-1$) all act as zero on $2\bj$
and $2\bj - 4\e_{n-1}$.  The operator $H_{n-1}$ (corresponding to $\a_{n-1}$) 
acts as zero on $2\bj$ and acts as the identity on $2\bj - 4\e_{n-1}$.
The operator $H_n$ (corresponding to $\a'_n$) acts as the identity on
$2\bj$ and as zero on $2\bj - 4\e_{n-1}$.
The second assertion is then proved
by adapting the corresponding argument in Proposition \secd.3 (iii).  
\qed\enddemo

\head \secf. The hyperoctahedron \endhead

In \S\secf, we consider examples relating to the polytope known as the
the {\it hyperoctahedron} or {\it cross polytope}; in Coxeter's notation
it is denoted $\be_n$.  The set $\Psi$ defined in Lemma \secf.1 is
our standard set of coordinates for the $2n$ vertices of the hyperoctahedron.

We will show how to use the hyperoctahedron to construct the remaining
two types of minuscule representations, namely the natural representations
for Lie algebras of types $C_n$ and $D_n$.

\proclaim{Lemma \secf.1}
Let $n \geq 4$, let
$\e_0, \ldots, \e_{n-1} \in \real^n$ be the usual basis for
$\real^n$, and let $$
\Psi = \{ \pm 4 \e_i : 0 \leq i \leq n-1 \}
.$$  Let $\Delta_D$ be as in Lemma \sece.3.  Then $\Psi$ is a minuscule
system with respect to $\Delta_D$.
\endproclaim

\demo{Proof}
We check Definition \secb.1, treating each vector $\ba \in \Delta_D$ in turn.
Suppose first that $\ba = \a_i$ for some $1 \leq i \leq n-1$, and let $\bv
\in \Psi$.

Define $\e_j$ to be the unique basis element such that $\bv = \pm 4\e_j$.
If $j \not\in \{i-1, i\}$ then we have $c = c(\bv, \ba) = 0$ and neither
vector $\bv \pm \ba$ lies in $\Psi$, satisfying Definition \secb.1 (ii).
If $\bv \in \{4\e_{i-1}, -4\e_i\}$ then $\bv - \ba \in \Psi$, 
$2 \bv . \ba = 32 = \ba . \ba$, giving $c = 1$ as required.  
The other possibility is that
$\bv \in \{-4\e_{i-1}, 4\e_i\}$, in which case $\bv + \ba \in \Psi$,
$2 \bv . \ba = -32 = -\ba . \ba$, giving $c = -1$ as required.

Now suppose $\ba = \a'_n$.  In this case, 
if $j \not\in \{n-2, n-1\}$ then we have $c = c(\bv, \ba) = 0$ and neither
vector $\bv \pm \ba$ lies in $\Psi$.
If $\bv \in \{4\e_{n-2}, 4\e_{n-1}\}$ then $\bv - \ba \in \Psi$, 
$\bv + \ba \not\in \Psi$, and
$2 \bv . \ba = 32 = \ba . \ba$, giving $c = 1$ as required.  
The other possibility is that
$\bv \in \{-4\e_{n-2}, -4\e_{n-1}\}$, in which case $\bv + \ba \in \Psi$,
$\bv - \ba \not\in \Psi$, and
$2 \bv . \ba = -32 = -\ba . \ba$, giving $c = -1$ as required.

Finally, suppose $\ba = \a_0$.  In this case, 
if $j \not\in \{0, 1\}$ then we have $c = c(\bv, \ba) = 0$ and neither
vector $\bv \pm \ba$ lies in $\Psi$.
If $\bv \in \{-4\e_0, -4\e_1\}$ then $\bv - \ba \in \Psi$, 
$\bv + \ba \not\in \Psi$, and
$2 \bv . \ba = 32 = \ba . \ba$, giving $c = 1$ as required.  
The other possibility is that
$\bv \in \{4\e_{i-1}, 4\e_i\}$, in which case $\bv + \ba \in \Psi$,
$\bv - \ba \not\in \Psi$, and
$2 \bv . \ba = -32 = -\ba . \ba$, giving $c = -1$ as required.
\qed\enddemo

\proclaim{Proposition \secf.2}
Maintain the notation of Lemma \secf.1.
\item{\rm (i)}{The $2n$-dimensional $\complex$-vector space $V_\Psi$ 
has the structure of a $\fg$-module, where $\fg$ is the derived affine 
Kac--Moody algebra of type $D_n^{(1)}$.}
\item{\rm (ii)}{Let $\Psi' = \Psi$ and $\Delta' = \Delta \backslash \{\a_0\}$.
Then $\Psi'$ is a minuscule system with respect to the simple system
$\Delta'$, and $V_{\Psi'}$ is a module for the simple Lie algebra $\fd_n$ 
over $\complex$ of 
type $D_n$.  It is an irreducible module with highest weight vector $4 \e_0$
and lowest weight vector $-4 \e_0$,
and affords the natural representation of $\fd_n$.}
\endproclaim

\demo{Proof}
Using Lemma \secf.1 in place of Lemma \secd.2, the proof
of (i) follows the same argument as the proof of Proposition \secd.3 (i).

The first assertion of (ii) follows by using Lemma \secf.1 and
copying the argument of Proposition \secd.3 (ii).  
The operators $H_i$ (for $1 < i \leq n$, where $H_n$ corresponds to
$\a'_n$) all act as zero on $4 \e_0$.  The operator $H_1$ acts as the identity
on $4 \e_0$.
The second assertion is then proved
by adapting the corresponding argument in Proposition \secd.3 (ii).  
\qed\enddemo

\proclaim{Lemma \secf.3}
Let $n \geq 2$, let $\Psi$ be as in Lemma \secf.1, and let $$
\Delta_C = \{\a_1, \ldots, \a_{n-1}\} \cup \{\a''_0, \a''_n\}
,$$ where $\a_i$ is as in Lemma \sece.1 for $1 \leq i \leq n-1$, 
$\a''_0 = -8 \e_0$ and $\a''_n = 8 \e_{n-1}$.
Then $\Psi$ is a minuscule system with respect to $\Delta_C$.
\endproclaim

\demo{Proof}
We check Definition \secb.1, treating each vector $\ba \in \Delta_C$ in turn.
The only cases not already covered by Lemma \secf.1 are the cases
where $\ba \in \{\a''_0, \a''_n\}$.  Let $\bv \in \Psi$, and
define $\e_j$ to be the unique basis element such that $\bv = \pm 4\e_j$.

Suppose that $\ba = \a''_n$.
If $j \ne n-1$ then we have $c = c(\bv, \ba) = 0$ and neither
vector $\bv \pm \ba$ lies in $\Psi$, satisfying Definition \secb.1.
If $\bv = \pm 4 \e_{n-1}$ then $\bv \mp \ba \in \Psi$, 
$\bv \pm \ba \not\in \Psi$,
and $2 \bv . \ba = \pm 64 = \pm \ba . \ba$, giving $c = \pm 1$ as required.  

The other possibility is that $\ba = \a''_0$.
If $j \ne 0$ then we have $c = c(\bv, \ba) = 0$ and neither
vector $\bv \pm \ba$ lies in $\Psi$, satisfying Definition \secb.1.
If $\bv = \pm 4 \e_0$ then $\bv \pm \ba \in \Psi$, $\bv \mp \ba \not\in \Psi$,
and $2 \bv . \ba = \mp 64 = \mp \ba . \ba$, giving $c = \mp 1$ 
and completing the proof.
\qed\enddemo

\proclaim{Proposition \secf.4}
Maintain the notation of Lemma \secf.3.
\item{\rm (i)}{The $2n$-dimensional $\complex$-vector space $V_\Psi$ 
has the structure of a $\fg$-module, where $\fg$ is the derived affine 
Kac--Moody algebra of type $C_n^{(1)}$.}
\item{\rm (ii)}{Let $\Psi' = \Psi$ and $\Delta' = \Delta \backslash 
\{\a''_0\}$. 
Then $\Psi'$ is a minuscule system with respect to the simple system
$\Delta'$, and $V_{\Psi'}$ is a module for the simple Lie algebra $\fc_n$ 
over $\complex$ of 
type $C_n$.  It is an irreducible module with highest weight vector $4 \e_0$
and lowest weight vector $-4 \e_0$,
and affords the natural representation of $\fc_n$.}
\endproclaim

\demo{Proof}
The proof is the same as the proof of Proposition \secf.2, using 
Lemma \secf.3 in place of Lemma \secf.1.
\qed\enddemo

\head \secg. Lines on Del Pezzo surfaces \endhead

In \S\secg, we revisit the examples of \S\secd\ involving the exceptional
Lie algebras $\fe_6$ and $\fe_7$.  We will highlight the close link
between the representation theory and the combinatorial algebraic
geometry associated with configurations of lines on Del Pezzo surfaces.
For more details on the latter, the reader is referred to \cite{{\bf 9}, \S V.4}.

\proclaim{Lemma \secg.1}
Let $\Psi$ and $\Delta$ be as in Proposition \secd.3, and let $K_0$ and $K_7$
be as in Definition \secd.1.
The action of the generators $\{s_\ba: \ba \in \Delta\}$ of the Weyl group 
$W = W_{\Psi, \Delta}$ on $\Psi$ are as follows.  If $\ba = \a_i$ with 
$0 \leq i \leq 6$, then $$s_\ba(\pm \bv_{j, k}) = s_\ba(\pm 
\bv_{s_i(j), s_i(k)}),$$
where $s_i$ is the simple transposition $(i, i+1)$.  We have 
$s_{\a_7}(\pm \bv_{i, j}) = \pm \bv_{i, j}$ unless $\{i, j\} \subset K_k$
for some $k \in \{0, 7\}$, in which case we have $$
s_{\a_7}(\pm \bv_{\{i, j\}}) = \mp \bv_{K_k \backslash \{i, j\}}
.$$  The action of $W$ on $\Psi$ is transitive.
\qed\endproclaim

\demo{Proof}
The formulae for the action of the $s_\ba$ are obtained by a routine case 
by case check.

The action of the $s_{\a_i}$ for $0 \leq i \leq 6$ makes it clear that
the vectors 
$\{+ \bv_{i, j} : 0 \leq i < j \leq 7\}$ 
are $W$-conjugate to each other, as are the vectors 
$\{- \bv_{i, j} : 0 \leq i < j \leq 7\}$.  The fact that $s_{\a_7}(+\bv_{0, 1})
= -\bv_{2, 3}$ completes the proof.
\qed\enddemo

\remark{Remark \secg.2}
The transformations induced by $s_{\a_7}$ described in the preceding proof
are sometimes known as {\it bifid transformations} (see Example 3 of 
\cite{{\bf 13}, \S4}).
\endremark

\proclaim{Lemma \secg.3}
Let $\Psi$ and $\Delta$ be as in Proposition \secd.3.
The diagonal action of the Weyl group 
$W$ on $\Psi \times \Psi$ has four orbits, each of which consists of a set $$
\{(\bv_1, \bv_2) : \bv_1, \bv_2 \in \Psi \text{ and } |\bv_1 - \bv_2| = D\}
$$ for some fixed number $D$.  More explicitly, the orbits are as follows:
\item{\rm (i)}{$\{(\bv, \bv) : \bv \in \Psi\}$, corresponding to $D = 0$;}
\item{\rm (ii)}{$\{(\pm\bv_{i, j}, \pm\bv_{i, k}) : |\{i, j, k\}| = 3\}
\cup 
\{(\pm\bv_{i, j}, \mp\bv_{k, l}) : |\{i, j, k, l\}| = 4\}$, corresponding to
$D = \sqrt{32}$,}
\item{\rm (iii)}{$\{(\pm\bv_{i, j}, \mp\bv_{i, k}) : |\{i, j, k\}| = 3\}
\cup 
\{(\pm\bv_{i, j}, \pm\bv_{k, l}) : |\{i, j, k, l\}| = 4\}$, corresponding to
$D = \sqrt{64}$,}
\item{\rm (iv)}{$\{(\bv, -\bv) : \bv \in \Psi\}$, corresponding to 
$D = \sqrt{96}$.}
\endproclaim

\demo{Proof}
The assertions about $D$ are easy to check.  This other assertions, which 
are also not difficult to prove, are a restatement of \cite{{\bf 4}, (4.1)}.
\qed\enddemo

\proclaim{Proposition \secg.4}
The $56$ elements of $\Psi$ are in natural bijection with the $56$ lines of 
the Del Pezzo surface of degree $2$; more precisely, 
if $\bv_1, \bv_2 \in \Psi$ are distinct points
with $|\bv_1 - \bv_2| = \sqrt{32 D}$, then $D - 1$ is the intersection number 
of the lines corresponding to $\bv_1$ and $\bv_2$.  In particular, pairs of
points at distance $\sqrt{32}$ correspond to skew lines on the Del Pezzo 
surface.

The $27$ elements of $\Psi(\bn, 8)$ (defined in Proposition \secd.3 (iii)) 
are $$
\{\bv_{0, i} : 1 \leq i \leq 6\} 
\cup
\{-\bv_{i, j} : 1 \leq i < j \leq 6\}
\cup
\{\bv_{i, 7} : 1 \leq i \leq 6\}
.$$  These are in natural bijection with the $27$ lines of the Del
Pezzo surface of degree $1$: in Hartshorne's notation \cite{{\bf 9}, Theorem 
V.4.9}, we identify $E_i$ with $\bv_{0, i}$, $F_{ij}$ with $-\bv_{i, j}$ and
$G_i$ with $\bv_{i, 7}$.  
The intersection number is defined as in the case of the $56$ lines.
\endproclaim

\demo{Proof}
The assertions about the Del Pezzo surface of degree $2$ are proved in 
\cite{{\bf 6}, p28}, where it is shown that $$
|\bv_1 - \bv_2|^2 = d^2 (x + 1)
,$$ where $\bv_1$ and $\bv_2$ are two distinct points of $\Psi$, $d$ is the 
minimal nontrivial distance between two points, and $x$ is the intersection 
number of the pair of lines corresponding to $\bv_1$ and $\bv_2$.  (The
precise link with the polytopes $2_{21}$ and $3_{21}$ is given on
\cite{{\bf 6}, p33}.)

It is easily checked that the $27$ elements of $\Psi(\bn, 8)$ are as listed.
By the result mentioned above, the only possible intersection numbers 
for two distinct lines on the Del Pezzo surface of degree $3$ are $0$ (meaning
the lines are skew) and $1$ (meaning the lines are incident).  Since no
two elements of $\Psi(\bn, 8)$ are at distance $\sqrt{96}$, it remains
to check that two distinct points of $\Psi(\bn, 8)$ are at distance 
$\sqrt{32}$ if and only if the corresponding lines are skew, and this 
follows from the rules given in \cite{{\bf 9}, Remark V.4.10.1}.
\qed\enddemo

Note that, because $\Psi(\bn, -8) = -\Psi(\bn, 8)$, the two $27$-dimensional
representations of $\fe_6$ are interchangeable in this context.

The next result explains how to recover the root system of type $E_7$ from
the set $E$ of directed edges of the polytope $3_{21}$.

\proclaim{Proposition \secg.5}
Maintain the notation of Lemma \secg.3.  Let $$
E = \{ (\bv_1, \bv_2) : \bv_1, \bv_2 \in \Psi \text{ and } |\bv_1 - \bv_2| =
\sqrt{32} \}
$$ has size $1512$.  The vectors 
$E' = \{\bv_1 - \bv_2 : (\bv_1, \bv_2) \in E\}$
form a root system of type $E_7$, and each of the $126$ roots occurs
with multiplicity $12$ in $E$.
\endproclaim

\demo{Proof}
If $\bv_1 = \bv_{0, 1}$ then one checks directly that there are $27$
vectors $\bv_2$ such that $(\bv_1, \bv_2) \in E$.  The fact (Lemma \secg.1)
that $W$ acts transitively on $\Psi$ implies by Lemma \secg.1 that $E$ has 
size $|\Psi| \times 27 = 1512$.

For the second assertion, note that $\bv_{0, 1} - \bv_{0, 2} = \a_1$.
The (additive) action of $W$ on $\Psi$ induces an action on $E'$, and
Lemma \secg.3 (ii) shows that $W$ acts transitively on $E'$.
Note that if $A_{ij} = A_{ji} = -1$, then $s_i s_j (\a_i) = \a_j$.
This implies that all the roots $\a_i$ are conjugate under the action of the 
Weyl group, and then \cite{{\bf 11}, \S5.1} shows that
the orbit $W . \a_1$ consists precisely of the
root system of type $E_7$.  By transitivity of the action of $W$ on the
root system, each root in $E'$
occurs with the same multiplicity, and by \cite{{\bf 2}, Appendix} there are 
$126$ roots of type $E_7$.  Since $1512/126 = 12$, the proof is completed.
\qed\enddemo

\proclaim{Proposition \secg.6}
Let $\Psi'$ and $\Delta'$ be as in Proposition \secd.3, and let 
$V_{\Psi'}$ be the
corresponding $56$-dimensional representation of the Lie algebra $\fe_7$.
If $\bv \in \Psi'$ and $x \in \fe_7$, then we have $$
x . \bv = \sum_{\bu \in \Psi''} \l_\bu \bu
,$$ where $\Psi'' = \{\bv\} \cup \{\bu : |\bv - \bu| = \sqrt{32}\}$.
In other words, if $\l_\bu \ne 0$, then either $\bu = \bv$ or the lines
on the Del Pezzo surface of degree $2$ corresponding to $\bu$ and $\bv$ are
skew.

A similar result holds for either of the $27$-dimensional representations of
$\fe_6$ and the Del Pezzo surface of degree $3$.
\endproclaim

\demo{Proof}
By \cite{{\bf 2}, \S4.1}, we have $$
\fe_7 = \fh \oplus \bigoplus_{\a \in \Phi} \fg_\a
,$$ where $\fh$ is a $7$-dimensional Cartan subalgebra, $\Phi$ is the root 
system for $\fe_7$, and the subspaces $\fg_\a$ are one-dimensional.  We 
identify $\fe_7$ with the algebra of operators on the $56$-dimensional
module $V$ as described in Proposition \secd.3 (ii).  

With these identifications,
if $\a = \a_i$ for $i \ne 0$, then $\fg_\a$ (respectively,
$\fg_{-\a}$ is spanned by the Lie algebra element 
$E_i = E_{\a_i}$ (respectively, $F_i = F_{\a_i}$).  The Cartan subalgebra
$\fh$ has as a basis the operators $H_i = H_{\a_i}$ for $1 \leq i \leq 7$.

It is possible to extend this to a basis for $\fe_7$ in which
(a) the subspace $\fg_\a$ for $\a$ a positive root is spanned by a vector of
the form $$
[\cdots [[E_{i_1}, E_{i_2}] E_{i_3}] \cdots E_{i_m}]
$$ where $\a = \sum_{j = 1}^m \a_{i_j}$ and 
(b) the subspace $\fg_\a$ for $\a$ a negative root is spanned by a vector of
the form $$
[\cdots [[F_{i_1}, F_{i_2}] F_{i_3}] \cdots F_{i_m}]
$$ where $-\a = \sum_{j = 1}^m \a_{i_j}$.
(See \cite{{\bf 7}, Proposition 5.4 (ii), (iv)} or \cite{{\bf 11}, (7.8.5)} for more 
details.)  

It follows that if $b_\bv$ is a basis element of $V$, $\a \in \Phi$
and $g_\a \in \fg_\a$,
then $
g_\a . b_\bv = \l b_{\bv + \a}
$ for some scalar $\l$ (meaning that $g_\a . b_\bv = 0$ if 
$\bv + \a \not\in \Psi$).  If $\l \ne 0$, then Proposition \secg.5 shows
that the distance from $\bv$ to $\bv + \a$ is $\sqrt{32}$.  By Proposition
\secg.4, we see that $\bv$ and $\bv + \a$ correspond to skew lines.  It 
follows easily from the definition of the $H_i$ that if $h \in H$ then $
h . b_{\bv} = \l b_{\bv}
$ for some scalar $\l$.  

Combining these observations proves the assertions about the $56$-dimensional
representation.  The argument can be easily adapted to work for the 
$27$-dimensional representation, because the root system of type $E_6$
embeds naturally into the root system of type $E_7$.
\qed\enddemo

\head \sech. Concluding remarks \endhead

In the various constructions presented above for irreducible modules for
simple Lie algebras, we did not provide self-contained proofs that the
modules constructed were irreducible.  However, this was done only to save
space, and it is not hard to give an elementary field-independent proof 
that these modules are irreducible.

One application of the polytope approach to minuscule representations
is that one can describe the crystal graph of each of the irreducible
modules that arises from the construction directly in terms of the polytope.
To do this, one starts with the vertices of $\Psi$, and for each element
$\ba \in \Delta$ corresponding to a simple root of the simple Lie algebra, one
connects two vertices $\bv_1$ and $\bv_2$ of $\Psi$ by an edge labelled $\ba$
if $\bv_1 - \bv_2 = \ba$.  It is not hard to show that this produces a
realization of the crystal graph, with the extra property that two edges are 
parallel if and only if they have the same label.

In the cases where the pair $(\Psi, \Delta)$ corresponds to a representation
of a simple Lie algebra, the elements of $\Psi$ may be partially ordered
by stipulating that $\bv_1 \leq \bv_2$ if $\bv_2 - \bv_1$ is a positive
linear combination of elements of $\Delta$; this corresponds to the usual
partial order on the weights of a representation.  The resulting partial order
on $\Psi$ makes $\Psi$ into a distributive lattice under the operations of
greatest lower bound and least upper bound.  (This is not a priori obvious,
but follows, for example, from the full heaps approach; see \cite{{\bf 7}, 
Corollary 2.2}.)
It would be interesting to know whether there is an easy way to define
the meet and join operations directly from the data $(\Psi, \Delta)$.

It may be tempting to think that one can describe a basis for each of the 
simple Lie algebras described in this paper by including operators
$E_\ba$ and $F_\ba$ for {\it every} positive root $\ba$.  However, such an 
algebra of operators would not be closed under the Lie bracket (except in
trivial cases) and what is needed instead is to modify the definition of
these new operators to introduce sign changes in certain places.  We do not
know if there is an easy way to keep track of these signs using the formalism
developed in this paper, although there is a good way to do it in the full
heaps approach, using the notion of the ``parity'' of a heap; see \cite{{\bf 7},
Definition 4.3, Definition 6.3} for details.

\leftheadtext{} \rightheadtext{}
%\vfill\eject

\end